\documentclass{llncs}
\usepackage[utf8]{inputenc}
\usepackage{graphicx}      % include this line if your document contains figures
\usepackage{amsmath} 
\usepackage{todonotes}
\usepackage{amssymb}

\usepackage[numbers]{natbib}

\newtheorem{thm}{Theorem}

\newcommand{\vol}{\mu_0}

\newcommand{\Diff}{\mathrm{Diff}}

\newcommand{\Diffvol}{{\operatorname{SDiff}}}

\newcommand{\Prob}{\mathrm{Prob}}

\newcommand{\GI}{G^{I}}
\newcommand{\ph}{\varphi}
\newcommand{\id}{\operatorname{id}}
\newcommand{\grad}{\nabla}
\begin{document}

\title{Diffeomorphic random sampling using optimal information transport}

\author{Martin Bauer\inst{1} \and Sarang Joshi\inst{2} \and Klas Modin\inst{3}} 
% \authorrunning{Ivar Ekeland et al.} % abbreviated author list (for running head) % %%%% list of authors for the TOC (use if author list has to be modified) 
\tocauthor{Martin Bauer, Sarang Joshi, Klas Modin} 
 \institute{Department of Mathematics,
Florida State University\\ \email{bauer@math.fsu.edu}
 \and Department of Bioengineering, Scientific Computing and Imaging Institute, University of Utah\\ \email{sjoshi@sci.utah.edu} \and 
 Department of Mathematical Sciences, Chalmers University of Technology and the University of Gothenburg \\ \email{klas.modin@chalmers.se}
}

\maketitle              % typeset the title of the contribution

\begin{abstract}
In this article we explore an algorithm for diffeomorphic random sampling of nonuniform probability distributions on Riemannian manifolds.
The algorithm is based on \emph{optimal information transport} (OIT)---an analogue of optimal mass transport (OMT).
% Contrary to OMT, based on the Wasserstein distance, the underlying geometry in OIT is based on the Fisher--Rao metric of information geometry. 
Our framework uses the deep geometric connections between the Fisher-Rao metric on the space of probability densities and the right-invariant \emph{information metric} on the group of diffeomorphisms.
The resulting sampling algorithm is a promising alternative to OMT, in particular as our formulation is semi-explicit, free of the nonlinear Monge--Ampere equation.
Compared to Markov Chain Monte Carlo methods, we expect our algorithm to stand up well when a large number of samples from a low dimensional nonuniform distribution is needed.

% might be more efficient for low dimensional distributions where 

% This allows us to develop methods that are inherently more efficient then optimal transport based algorithms for diffeomorphic density registration.

\textbf{Keywords:} density matching, information geometry, Fisher--Rao metric, optimal transport, image registration, diffeomorphism groups, random sampling

\textbf{MSC2010:} 58E50, 49Q10, 58E10
 \end{abstract}

\section{Introduction}
We construct algorithms for random sampling, addressing the following problem.
\begin{problem}\label{pro:original}
	Let $\mu$ be a probability distribution on a manifold $M$. 
	Generate $N$ random samples from $\mu$. 
\end{problem}
The classic approach to sample from a probability distribution on a higher dimensional space is to use Markov Chain Monte Carlo (MCMC) methods, for example the Metropolis--Hastings algorithm~\cite{Ha1970}. 
An alternative idea is to use diffeomorphic density matching between the density $\mu$ and a standard density $\vol$ from which samples can be drawn easily.
Standard samples are then transformed  by the diffeomorphism to generate non-uniform samples.
 % which allows then to transform samples from the uniform distribution to generate a sample of the distribution $\mu$. 
In Bayesian inference, for example, the distribution $\mu$ would be the posterior distribution and $\vol$ would be the prior distribution.
In case the prior itself is hard to sample from the uniform distribution can be used. 
For $M$ being a subset of the real line, the standard approach is to use the cumulative distribution function to define the diffeomorphic transformation. 
If, however, the dimension of $M$ is greater then  one there is no obvious change of variables to transform the samples to the distribution of the prior.
We are thus led to the following matching problem.

\begin{problem}\label{prob:diff_dens}
	Given a probability distribution $\mu$ on $M$, find a diffeomorpism $\varphi$ such that
	\begin{equation*}
		\varphi_*\vol = \mu.
	\end{equation*}
	Here, $\vol$ denotes a standard distribution on $M$ from which samples can be drawn, and $\varphi_*$ is the the push-forward of $\varphi$ acting on densities, \emph{i.e.},
	\begin{equation*}
		\varphi_*\vol = |D\varphi| \vol\circ\varphi,
	\end{equation*}
	where $|D\varphi|$ is the Jacobian determinant.
\end{problem}
A benefit of transport-based methods over traditional MCMC methods is cheap computation of additional samples;
it amounts to drawing uniform samples and then evaluating the transformation.
On the other hand, transport-based methods scale poorly with increasing dimensionality of $M$, contrary to MCMC.

The action of the diffeomorphism group on the space of smooth probability densities is transitive (Moser's lemma~\cite{Mo1965}), so existence of a solution to Problem~\ref{prob:diff_dens} is guaranteed. 
%% Can we assume the above statement is obvious?
However, if the dimension of $M$ is greater then one, there is an infinite-dimensional space of solutions.
Thus, one needs to select a specific diffeomorphism within the set of all solutions. 
Moselhy and Marzouk \cite{MoMa2012} and Reich \cite{Re2013} proposed to use optimal mass transport (OMT) to construct the desired diffeomorphism $\varphi$, thereby enforcing $\varphi = \nabla c$ for some convex function~$c$. 
The OMT approach implies solving, in one form or another, the heavily non-linear Monge--Ampere equation for $c$.
A survey of the OMT approach to random sampling is given by Marzouk \emph{et.~al.}~\cite{MaMoPaSp2016}.

In this article we pursue an alternative approach for diffeomorphic based random sampling, replacing OMT by \emph{optimal information transport} (OIT), which is diffeomorphic transport based on the Fisher--Rao geometry~\cite{Mo2015}.
Building on deep geometric connections between the Fisher--Rao metric on the space of probability densities  and the right-invariant \emph{information metric} on the group of diffeomorphisms \cite{KLMP2013,Mo2015}, we developed in \cite{BJM2015} an efficient numerical method for density matching.
The efficiency stems from a solution formula for $\varphi$ that is explicit up to inversion of the Laplace operator, thus avoiding the solution of nonlinear PDE such as Monge--Ampere.
In this paper we explore this method for random sampling (the initial motivation in \cite{BJM2015} is medical imaging, although other applications, including random sampling, are also suggested).
The resulting algorithm is implemented in a short MATLAB code, available under MIT license at \texttt{\url{https://github.com/kmodin/oit-random}}.
% The formula is then the basis of an efficient numerical algorithm, as developed in .
% \todo[inline,backgroundcolor=blue!20!white,
% bordercolor=red]{We should mention more strongly that this is really based on our previous article. }

% , yielding and thus the obtained algorithms are by magnitudes faster then those based on OMT.

%
\section{Density Transport Problems}
Let $M$ be an $d$--dimensional orientable, compact manifold equipped with a Riemannian metric $g=\langle.,.\rangle$.
The volume density induced by $g$ is denoted $\vol$ and without loss of generality we assume that the total volume of $M$ with respect to $\vol$ is one, \emph{i.e.}, 
$\int_M \vol=1$.
Furthermore, the space of smooth probability densities on $M$ is given by
\begin{equation}
	\operatorname{Prob}(M) = \{ \mu\in\Omega^{d}(M)\mid \mu>0, \quad \int_M\mu = 1 \},
\end{equation}
where $\Omega^{d}(M)$ denotes the space of smooth $d$-forms. 
The group of smooth diffeomorphisms $\operatorname{Diff}(M)$ acts on the space of probability densities via push-forward:
\begin{align}
\Diff(M)\times \Prob(M)&\mapsto \Prob(M)\\
(\varphi,\mu)&\rightarrow \varphi_*\mu\;.
\end{align}
By  a result of Moser  \cite{Mo1965} this action is transitive. 

We introduce the subgroup of volume preserving diffeomorphisms
\begin{equation}
	\Diffvol(M) = \{\varphi\in\Diff(M)\mid \varphi_*\vol = \vol \}\;.
\end{equation}
Note that  $\Diffvol(M)$ is the  isotropy group of $\vol$ with respect to the action of $\Diff(M)$.
The spaces $\operatorname{Prob}(M)$, $\Diff(M)$, and $\Diffvol(M)$ all have the structure of smooth, infinite dimensional Fréchet manifold.
Furthermore, $\Diff(M)$ and $\Diffvol(M)$ are infinite dimensional Fréchet Lie groups.
A careful treatment of these Fréchet topologies can be found in the work by Hamilton \cite{Ha1982}.
 
In the following we will focus our attention on the diffeomorphic density matching problem (Problem~\ref{prob:diff_dens}). 
A common approach to overcome the non-uniqueness in the solution is to add a regularization term to the problem.
That is, to search for a minimum energy solution that has the required matching property, for some energy functional $E$ on the diffeomorphism group.  
Following ideas from mathematical shape analysis \cite{MiTrYo2002} it is a natural approach to define this energy functional using the geodesic distance function $\operatorname{dist}$ of a Riemannian metric on the diffeomorphism group. 
Then the regularized diffeomorphic matching problem can be written as follows.
\begin{problem}\label{prob:oit}
	Given a probability density $\mu\in \operatorname{Prob}(M)$ we want to find the diffeomorphism $\varphi \in \operatorname{Diff}(M)$ that minimizes the energy functional
	\begin{align}
		E(\varphi)= \operatorname{dist}^{2}(\operatorname{id}, \varphi)
	\end{align}
	over all diffeomorphisms $\varphi$ with $\varphi_*\vol=\mu$.
\end{problem}
The free variable in the above matching problem is the choice of Riemannian metric---thus distance function---on the group of diffeomorphisms.
Although not formulated as here, Moselhy and Marzouk~\cite{MoMa2012} proposed to use the $L^2$ metric on $\Diff(M)$
	\begin{align}
		G_{\varphi}(u\circ\varphi,v\circ\varphi)= \int_M \langle u\circ\varphi,v\circ\varphi\rangle\; \vol\;
	\end{align}
for $u\circ\varphi,v\circ\varphi \in T_{\varphi}\Diff(M)$. 
This corresponds to distance-squared \emph{optimal mass transport} (OMT), which induces the Wasserstein $L^2$ distance on $\Prob(M)$, see, for example,~\cite{Ot2001,KhWe2009,Vi2009}. 

In this article we use the right-invariant $H^1$-type metric
\begin{equation}\label{eq:H1dot}
	\GI_{\varphi}(u\circ\varphi,v\circ\varphi )= -\int_M \langle\Delta u,v\rangle\vol + \lambda\sum_{i=1}^{k}\int_{M}\langle u,\xi_i\rangle\vol\, \int_{M}\langle v,\xi_i\rangle\vol,
	% \\ u\coloneqq U\circ\varphi^{-1},v\coloneqq V\circ\varphi^{-1} ,
\end{equation}
where $\lambda>0$, $\Delta$ is the Laplace--de~Rham operator lifted to vector fields, and $\xi_1,\ldots,\xi_k$ is an orthonormal basis of the harmonic 1-forms on~$M$.
Because of the Hodge decomposition theorem, $\GI$ is independent of the choice of orthonormal basis $\xi_1,\ldots,\xi_k$ for the harmonic vector fields.
% (see~\autoref{sub:notation} for details). 
This construction is related to the Fisher-Rao metric on the space of probability density \cite{Fr1991,BaBrMi2016}, which is predominant in the field of information geometry \cite{AmNa2000}. 
We call $\GI$ the \emph{information metric}. 
See \cite{KLMP2013,Mo2015,BJM2015} for more information on the underlying geometry.

The connection between the information metric and the Fisher-Rao metric allows us to construct almost explicit solutions formulas for Problem~\ref{prob:diff_dens} using the explicit formulas for the geodesics of the Fisher-Rao metric.	
\begin{thm}[\cite{Mo2015,BJM2015}]
Let $\mu\in \Prob(M)$ be a smooth probability density. 
The diffeomorphism $\varphi\in \Diff(M)$ minimizing 
 $\operatorname{dist}_{G^I}(\operatorname{id},\varphi)$ under the constraint $\varphi_*\vol=\mu$ is given by $\varphi(1)$, where $\varphi(t)$ is obtained as the solution to the problem
 \begin{equation}\label{eq:veq}
	\begin{split}
		\Delta f(t) &= \frac{\dot\mu(t)}{\mu(t)}\circ \ph(t), \\
		v(t) &= \nabla(f(t)), \\
		\frac{d}{d t}\ph(t)^{-1} &= v(t)\circ\ph(t)^{-1}, \quad \ph(0) = \id 
	\end{split}
\end{equation}
where $\mu(t)$ is the (unique) Fisher-Rao geodesic connecting $\vol$ and $\mu$
\begin{equation}\label{eq:fisher_rao_geodesics}
	\mu(t)= \left( 
			\frac{\sin\left((1-t)\theta\right)}{\sin\theta} + \frac{\sin\left(t\theta\right)}{\sin\theta}\sqrt{\frac{\mu}{\vol}}
		\right)^2 \vol,
		\quad \cos\theta = \int_M   \sqrt{\frac{\mu}{\vol}}\; \vol  \,.
\end{equation}
\end{thm}
% Note, that the above theorem shows -- by providing an explicit formula -- the existence and uniqueness of a solution to the regularized transport problem for an arbitrary target densities $\mu$. 
The algorithm for diffeomorphic random sampling, described in the following section, is directly based on solving the equations \eqref{eq:veq}. 

\section{Numerical Algorithm}
In this section we explain the algorithm for random sampling using optimal information transport.
It is a direct adaptation of \cite[Algorithm~1]{BJM2015}.

\noindent\rule{\textwidth}{0.4pt}\vspace{.15cm}\newline
\textbf{Algorithm 1} (OIT  based random sampling) \newline
\rule{\textwidth}{0.4pt}
Assume we have a numerical way to represent functions, vector fields, and diffeomorphisms on~$M$, and numerical methods for
\begin{itemize}
\item composing functions and vector fields with diffeomorphisms,
\item computing the gradient of functions, 
\item computing solutions to Poisson's equation on $M$,
\item sampling from the standard distribution $\vol$ on $M$, and
\item evaluating diffeomorphisms.
\end{itemize}
An OIT based algorithm for Problem~\ref{pro:original} is then given as follows:
\begin{enumerate}
\item Choose a step size $\varepsilon = 1/K$ for some positive integer~$K$ and calculate the Fisher-Rao geodesic $\mu(t)$ and its derivative $\dot \mu(t)$ at all time points $t_k=\frac{k}{K}$ using equation~\eqref{eq:fisher_rao_geodesics}.
	\item 
	Initialize $\varphi_0 = \id$. %, and $\varphi^{-1}_{0} = \id$.
	Set $k\leftarrow 0$.

	\item Compute $s_k = \frac{\dot \mu(t_k)}{\mu(t_k)}\circ\varphi_{k}$ and solve the Poisson equation
	\begin{equation}
		\Delta f_k= s_k.
	\end{equation}

	\item Compute the gradient vector field $v_k = \grad f_k$.

	\item Construct approximations $\psi_k$ to $\exp(-\varepsilon v_k)$, for example 
	\begin{equation}
		\psi_k = \id - \varepsilon v_k.		
	\end{equation}

	\item Update the diffeomorphism\footnote{If needed, one may also compute the inverse by $\varphi_{k+1}^{-1} = \varphi_k^{-1} + \varepsilon v\circ\varphi_k^{-1}$.}
	\begin{equation}
		% \varphi_{k+1} = \psi_k\circ\varphi_k, \quad 
		\varphi_{k+1} = \varphi_k\circ\psi_k .
	\end{equation}
	
	\item Set $k \leftarrow k+1$ and continue from step 3 unless $k=K$.
% \end{enumerate}
% The random samples $y_1,\ldots, y_N$ from $\mu$ are then obtained as follows:
% \begin{enumerate}
 % \setcounter{enumi}{7}
\item Draw $N$ random samples $x_1,\ldots x_N$ from the uniform distribution $\vol$.
\item Set $y_n=\varphi_K(x_n)$,  $n\in\{1,\ldots N\}$.
\end{enumerate}
\vspace{-.4cm}
\rule{\textwidth}{0.4pt}

\vspace{3ex}

The algorithm generates $N$ random samples $y_1,\ldots,y_N$ from the distribution~$\mu$.
One can save $\varphi_K$ and repeat 8-9 whenever additional samples are needed.

The computationally most intensive part of the algorithm is the solution of Poisson's equation at each time step.
Notice, however, that we do not need to solve nonlinear equations, such as Monge--Ampere, as is necessary in OMT.

% Note that the Poisson equation
% Note that in our density matching  algorithm one only has to solve a series of linear Poisson equation, compared to the OMT approach where solving the transport problem requires  to solve the non-linear Monge-Ampere equation.

\section{Example}

\begin{figure}[http]
	\includegraphics[height=0.3\textheight]{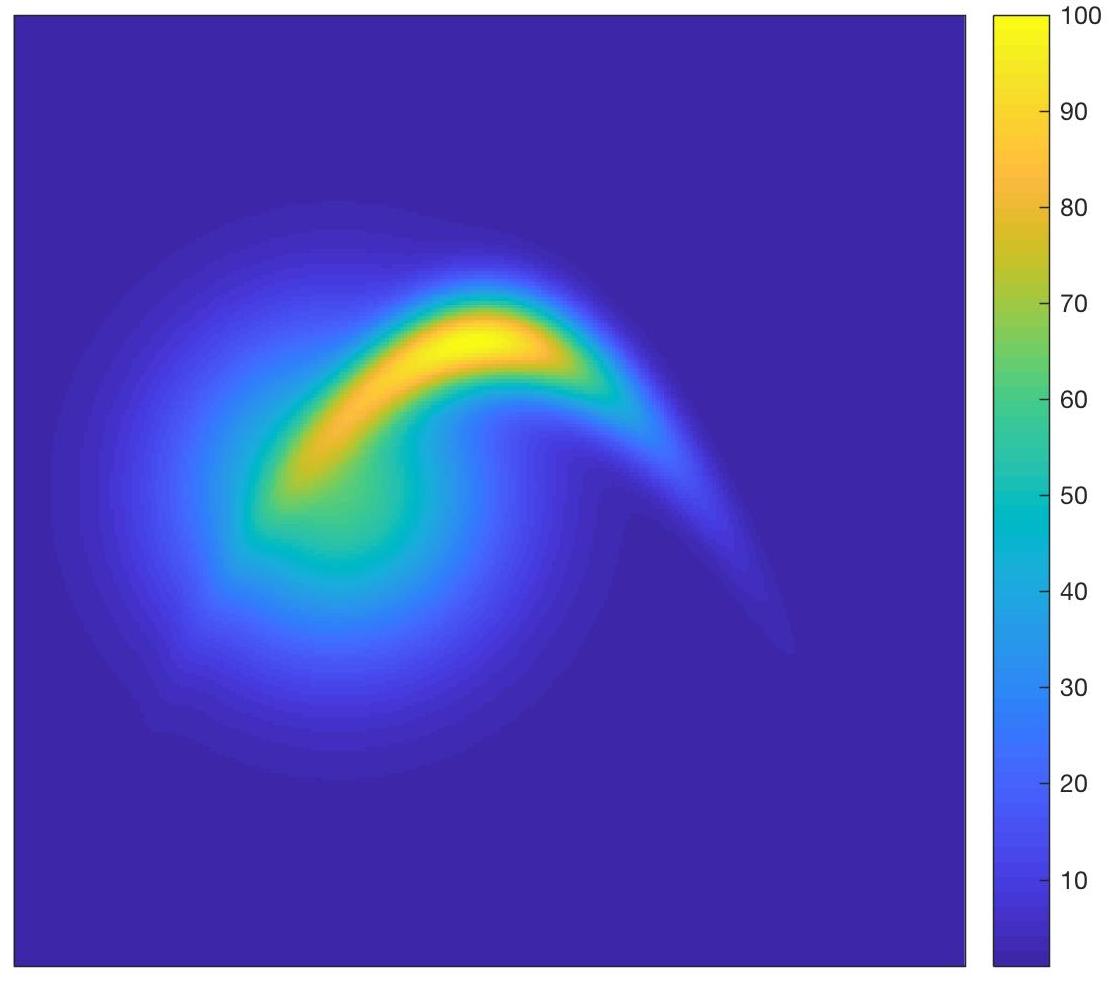}
	\includegraphics[height=0.3\textheight]{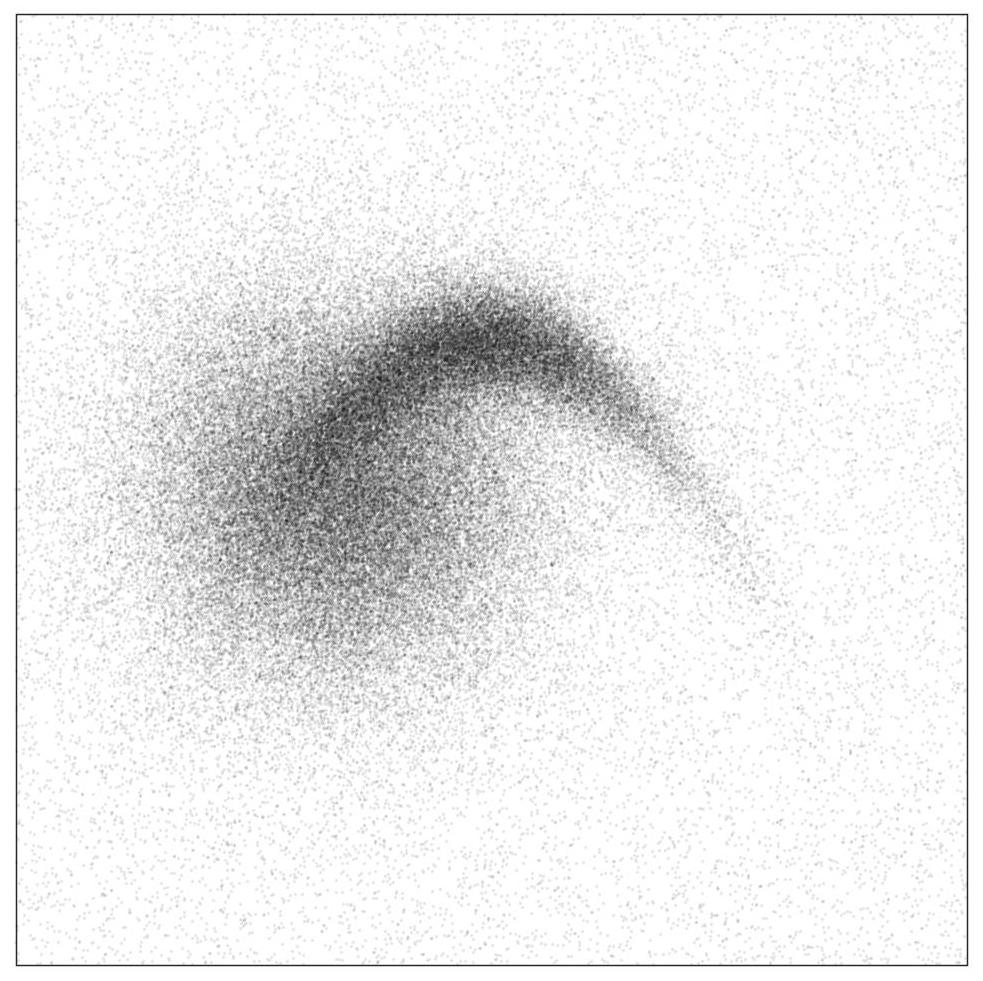}
	\caption{
	(left) The probability density $\mu$ of \eqref{distribution}.
	The maximal density ratio is 100.\qquad
	(right) $10^5$ samples from $\mu$ calculated using our OIT based random sampling algorithm.}
	\label{fig:example1_samples}
\end{figure}

In this example we consider $M=\mathbb{T}^2 \simeq (\mathbb{R}/2\pi\mathbb{Z})^2$ with distribution defined in Cartesian coordinates $x,y\in [-\pi,\pi)$ by
\begin{equation}\label{distribution}
	\mu \sim 3\exp(-x^2 - 10(y-x^2/2+1)^2) + 2\exp(-(x+1)^2 - y^2)+1/10,
\end{equation}
normalized so that the ratio between the maximum and mimimum of $\mu$ is 100.
The resulting density is depicted in Fig.~\ref{fig:example1_samples}~(left).

We draw $10^5$ samples from this distribution using a MATLAB implementation of our algorithm, available under MIT license at
\begin{center}
	\texttt{\url{https://github.com/kmodin/oit-random}}
\end{center}
% , with 100 time steps and spatial resolution of $256\times 256$.

\begin{figure}[http]
	\begin{center}
		\includegraphics[height=0.48\textheight]{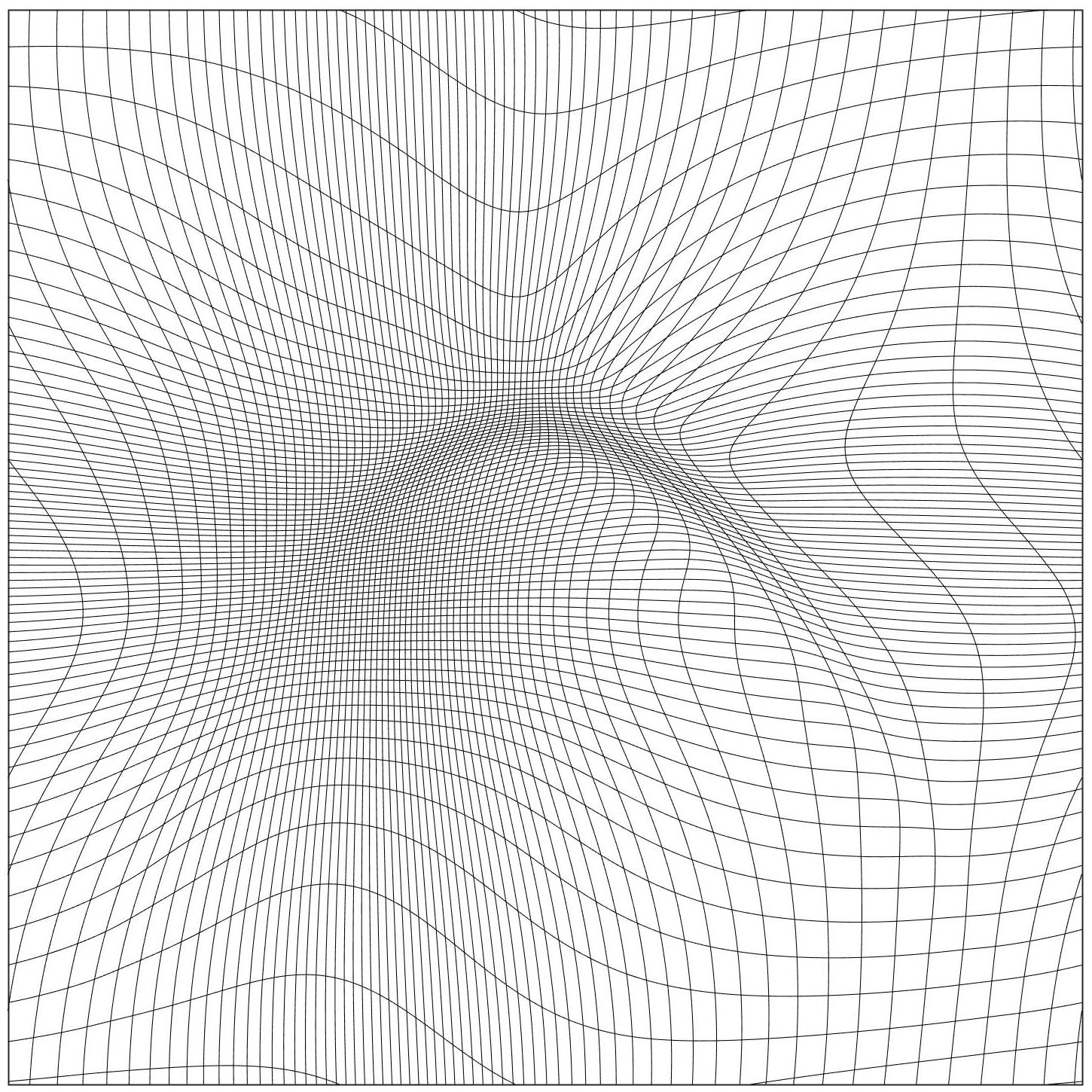}
	\end{center}
	\caption{
	The computed diffeomorphism $\varphi_K$ shown as a warp of the uniform $256\times 256$ mesh (every 4th mesh-line is shown). 
	Notice that the warp is periodic.
	It satisfies $\varphi_*\vol=\mu$ and solves Problem~\ref{prob:oit} by minimizing the information metric \eqref{eq:H1dot}.
	The ratio between the largest and smallest warped volumes is 100.
	}
	\label{fig:example1_phiinv}
\end{figure}

The implementation can be summarized as follows.
To  solve the lifting equations \eqref{eq:veq} we discretize the torus by a $256\times 256$ mesh and use the fast Fourier transform (FFT) to invert the Laplacian. 
We use 100 time steps.
The resulting diffeomorphism is shown as a mesh warp in Fig.~\ref{fig:example1_phiinv}.
We then draw $10^5$ uniform samples on $[-\pi,\pi]^2$ and apply the diffeomorphism on each sample (applying the diffeomorphism corresponds to interpolation on the warped mesh).
The resulting random samples are depicted in Fig.~\ref{fig:example1_samples}~(right).
To draw new samples is very efficient.
For example, another $10^7$ samples can be drawn in less than a second on a standard laptop.

%\subsection{Specific Examples of Transforming Manifolds: 2-Sphere $S^{2}$}
%Given an arbitrary probability distribution on the surface of the unit sphere we use the Spherical Harmonic expansion to solve analytically the spherical Laplacian. Analogous to the Fourier exponentials on the torus, the Spherical Harmonics are Eigen function of the Laplacian and Fast Spherical Harmonic transforms exists~\cite{MartinFastSphericalTransform} which gives us a very efficient algorithm for computing the flow. Using 
\vspace{-.5cm}
\section{Conclusions}
In this paper we explore random sampling based on the optimal information transport algorithm developed in \cite{BJM2015}.
% The resulting algorithm is well-posed on any compact Riemannian manifold, given that one has a way to numerically represent diffeomorphisms on the manifold.
Given the semi-explicit nature of the algorithm, we expect it to be an efficient competitor to existing methods, especially for drawing a large number of samples from a low dimensional manifold.
However, a detailed comparison with other methods, including MCMC methods, is outside the scope of this paper and left for future work.

We provide an example of a complicated distribution on the flat 2-torus.
The method is straighforward to extended to more elaborate manifolds, \emph{e.g.}, by using finite element methods for the efficient solution of Poisson's equation on manifolds.
For non-compact manifolds, most importantly $\mathbb{R}^n$, one might use standard techniques, such as Box--Muller, to first transform the required distribution to a compact domain.

\bibliographystyle{splncs03} % (plain, alpha, abbrv, acm, ieee, abbrvnat)
\bibliography{refs}
% \bibliography{/Users/moklas/Documents/Papers/References}

\end{document}